\documentclass[a4paper, 12pt]{article}

\usepackage[koi8-r,cp1251]{inputenc}
\usepackage[english,russian]{babel}
\usepackage{amsfonts,amssymb,amsmath, hyperref}
\usepackage[final]{epsfig}
\usepackage{graphicx}

\usepackage{amscd}

\begin{document}

\begin{center}
\Large{\textbf{On compact Hankel operators over compact Abelian groups}}
 \end{center}

 \begin{center}
 A. R. Mirotin
 \end{center}

 \begin{center}
 amirotin@yandex.ru
 \end{center}

\
 
\begin{center}
 КОМПАКТНЫЕ ОПЕРАТОРЫ ГАНКЕЛЯ НАД КОМПАКТНЫМИ АБЕЛЕВЫМИ ГРУППАМИ
\end{center}

\begin{center}
А. Р. Миротин
\end{center}

\

AMS SC 2010: Primary 43A17, secondary 47B35.

\

\textit{Key wards:} compact Abelian group, linearly ordered group, Hankel operator, compact operator, finite rank operator, Schatten--von Neumann class, singular number, invariant subspace.

\

We consider compact and connected Abelian group $G$ with
a linearly ordered dual.
  Based on the description of the structure of compact Hankel operators over $G$, generalizations of the classical Kronecker, Hartman, Peller and Adamyan - Arov - Krein theorems are obtained. A generalization of Burling's invariant subspace theorem is also established. Applications are given to Hankel operators over discrete groups

\

\begin{center}
 Введение
\end{center}

\

Как известно, пространство Харди $H^2$ на окружности является одной из наиболее удобных реализаций сепарабельного
 бесконечномерного гильбертова пространства. Операторы Ганкеля в этих пространствах, их различные реализации и аналоги активно
изучались  и нашли важные приложения, в том числе к теории функций, теории операторов, теории случайных процессов, теории
 рациональных приближений и теории управления
 (см., например, \cite{Pel} --- \cite{Nik2}). Изучались также многочисленные   обобщения этих операторов
(см. \cite{Nik1} --- \cite{YCG} и обзор в \cite[с. 195 -- 204]{Nik1}), и обобщения тесно связанных с ними операторов Тёплица
 и операторов Винера-Хопфа (см. \cite{Adukov}, \cite{SbMath},  \cite{EMRS} и библиографию там).

Если классические ганкелевы операторы действуют в сепарабельных гильбертовых пространствах, то   обобщенные операторы Ганкеля,
рассматриваемые в данной работе, определены в гильбертовых пространствах, которые не обязаны быть сепарабельными, в частности,
в пространствах   $H^2(G),$ где
   $G$ есть нетривиальная связная  компактная абелева группа с (дискретной) группой характеров $X.$ Известно \cite{Pont}, что  группа характеров $X$
связной (и только связной) компактной абелевой группы может быть наделена линейным порядком,
согласованным со структурой группы, и не имеет кручения (см., например, \cite{Rud}). При этом линейный порядок
в $X,$ вообще говоря, не единственен. Ниже мы предполагаем, что такой порядок выбран и фиксирован.
Это равносильно тому, что в группе  $X$ выделена   подполугруппа $X_+$ (положительный конус),  содержащая
единичный характер $1$
и такая, что $X_+\cap X_+^{-1}=\{1\}$ и $X=X_+\cup X_+^{-1}$. При этом
полугруппа $X_+$ индуцирует в $X$  линейный порядок, согласованный
со структурой группы,  по правилу $\xi\leq\chi:=\chi\xi^{-1}\in
X_+$. Ясно, что  $X_+=\{\chi\in X:\chi\geq 1\}$. Далее мы положим $X_-:=X\setminus X_+$
($= X_+^{-1}\setminus\{1\}$).

В приложениях в роли  $X$ часто выступают подгруппы аддитивной группы $\mathbb{R}^n,$
наделенные дискретной топологией, так что $G$ является боровской компактификацией группы  $X$ (см, например, \cite{EMRS}
 и цитированную там литературу). В частности, в качестве  $X$ можно взять группу $\mathbb{Z}^n,$ наделенную лексикографическим
 порядком. В этом случае $X$ имеет наименьший положительный элемент, а $G=\mathbb{T}^n$ ---  $n$-мерный тор.
 (Отметим, что описание всех линейных порядков на группе $\mathbb{Z}^n$ дано в  \cite{Teh}, \cite{Zajt}.)  Другими интересными
 примерами могут служить бесконечномерный тор $\mathbb{T}^\infty$ (см. \cite[пример 3]{SbMath}, где в этом случае указан линейный порядок на  $X,$ при котором существует наименьший положительный элемент), боровский компакт (группа характеров группы $\mathbb{R},$
  наделенной дискретной топологией и естественным порядком), а также группы характеров аддитивной группы $\mathbb{Q}$ рациональных чисел и мультипликативной группы $\mathbb{Q}_{>0}$ положительных рациональных чисел, наделенных дискретной топологией и естественным порядком. В последних  трех примерах $X$  не имеет наименьшего положительного элемента. (В случае группы $\mathbb{Q}_{>0}$ под положительными элементами понимаются рациональные числа, большие 1, а $G$ есть бесконечномерный тор $\mathbb{T}^\infty.$ Таким образом, в группе характеров группы $\mathbb{T}^\infty$ имеется как линейный порядок с наименьшим положительным элементом \cite[пример 3]{SbMath},  так и линейный порядок без оного.)

  Данная работа продолжает цикл работ автора и его учеников \cite{MD1} --- \cite{MD4}. Используя другой подход, чем в
  \cite{MD1} --- \cite{MD4}, мы  для операторов Ганкеля в пространствах $H^2$ на
  компактных абелевых группах доказываем обобщения классических теорем Кронекера и Хартмана
(дающих, соответственно, критерии конечномерности и компактности операторов Ганкеля в пространствах $H^2$ на  окружности),
 сняв ограничения на оператор, наложенные в этих работах. Кроме того, для рассматриваемых  операторов Ганкеля мы устанавливаем
   обобщения классических теорем Пеллера и Адамяна--Арова--Крейна, дающих, соответственно, критерий принадлежности оператора
    классам Шаттена--фон Неймана и выражение для  сингулярных чисел компактного ганкелева оператора.
Эти результаты опираются на редукцию компактных операторов Ганкеля над произвольными компактными абелевыми группами с
линейно упорядоченной группой характеров к операторам Ганкеля над группой окружности $\mathbb{T}.$  Показано, что для компактного ганкелева оператора на группе символ фактически факторизуется через первый положительный характер, что сводит общий случай к случаю окружности. Тем самым, выясняется исключительная роль группы $\mathbb{T}$  в этом круге вопросов. Получено также обобщение теоремы Бёрлинга об инвариантных подпространствах, дополняющее известный результат Хелсона--Лауденслегера.

Резюмируя основные результаты работы можно сказать, что такие фундаментальные факты классического анализа на группе окружности как теорема Бёрлинга об инвариантных подпространствах, наличие конечных произведений Бляшке,  теоремы Кронекера,  Хартмана,  Пеллера и Адамяна--Арова--Крейна сохраняются для абелевых компактных связных групп (с фиксированным линейным порядком на двойственной группе) тогда и только тогда, когда в двойственной группе имеется первый положительный элемент.

Работа завершается  приложениями полученных результатов к
операторам Ганкеля над дискретными абелевыми группами.

\

\begin{center}
\S 1. Операторы Ганкеля над компактными абелевыми группами
\end{center}

\

Если в группе  $X$  имеется наименьший положительный элемент $\chi_1$, то через  $X^i$ будем обозначать
 (бесконечную циклическую) подгруппу группы  $X$, порожденную этим элементом    (см.  \cite[теорема 2]{SbMath}).
  Следующая лемма доказана в \cite{MD4}.

{\bf  Лемма 1.1.} \textit{ Пусть в группе  $X$  имеется наименьший положительный элемент} $\chi_1.$

1) \textit{ Множество $X_+\setminus X^i$  есть идеал полугруппы   $X_+;$}

2)   \textit{для любого $k\in \mathbb{Z}$ справедливо равенство}
$$
\chi_1^k(X_+\setminus X^i)=X_+\setminus X^i.
$$

{\bf Определение 1.2.} {\it Пространство Харди} $H^p(G)\ (1\leq p\leq\infty)$ {\it над}
$G$ определяется следующим образом (см., например, \cite{Rud}):
$$
H^p(G)=\{f\in L^p(G):\widehat{f}(\chi)=0\ \forall\chi\in X_-\},
$$
где $\widehat{f}$ обозначает преобразование Фурье функции  $f\in L^1(G).$ Норма в $L^p(G)$ (и в $H^p(G)$)
будет обозначаться $\|\cdot\|_p.$

Обозначим через $H^2_-(G)$ ортогональное дополнение подпространства
$H^2(G)$ пространства $L^2(G).$ Тогда
$$
H^2_-(G)=\{f\in L^2(G):\widehat{f}(\chi)=0\ \forall\chi\in X_+\}.
$$

Множество $X_+$ является ортонормированным базисом
пространства $H^2(G),$  а  $X_-$ ---
ортонормированным базисом пространства $H^2_-(G)$ (гильбертова размерность этих пространств равна мощности множества $X$).
Через $P_+$ и $P_-$ мы будем обозначать ортопроекторы из  $L^2(G)$ на $H^2(G)$ и $H^2_-(G)$ соответственно.
 Нормированная мера Хаара группы $G$ будет обозначаться через $dx.$

Отметим, что группы $X,$ для которых $X_+$ содержит наименьший элемент, описаны в \cite[лемма 3.2]{Adukov}.

{\bf Определение 1.3.} Пусть $\varphi\in L^{2}(G).$ {\it
\textit{Оператором Ганкеля} (ганкелевым оператором) над  $G$ с символом $\varphi$} назовем оператор
$H_{\varphi}:H^2(G)\rightarrow H^2_-(G),$ определяемый  равенством
$$
H_{\varphi}=P_-M_{\varphi},
$$
где  $M_{\varphi}: f\mapsto\varphi f$ --- оператор умножения на $\varphi.$

\textbf{Лемма 1.4.} $H_{\varphi}=H_{P_-\varphi}$ \textit{ при} $\varphi\in L^{2}(G).$

Доказательство. Достаточно показать, что $H_{\psi}=0$ при $\psi\in H^{2}(G).$
 Но если $\psi, f\in H^{2}(G),$ то их разложения в ряд Фурье имеют вид
 $\psi=\sum_{\chi\in X_+}a_\chi \chi,$ $f=\sum_{\xi\in X_+}b_\xi \xi,$ где
$(a_\chi), (b_\xi)\in \ell^2(X_+).$
Поэтому
$$
\psi f=\sum_{\chi\in X_+}\sum_{\xi\in X_+}a_\chi b_\xi \chi  \xi,
$$
причем
$$
\sum_{\chi\in X_+}\sum_{\xi\in X_+}|a_\chi b_\xi|^2=\left(\sum_{\chi\in X_+}|a_\chi|^2\right)
\left(\sum_{\xi\in X_+} |b_\xi|^2\right)<\infty.
$$
Так как
$\chi\xi\in X_+$ в разложении $\psi f,$ то  $\psi f\in H^2(G),$ что и завершает доказательство леммы.$\Box$

В статье \cite[c. 141]{MD1} определены пространства $BMO(G)$ и показано, что если   в группе $X$ существует наименьший
положительный элемент,  то при $\varphi\in L^2(G)$ оператор $H_{\varphi}$ ограничен тогда и только
тогда, когда $P_-\varphi\in BMO(G)$ (и, в частности, он ограничен, если $\varphi\in L^\infty(G)$),
и это равносильно тому, что $H_\varphi=H_g$ для некоторой функции $g,$ непрерывной на $G$ (см. там же, следствие 2).

Для доказательства обобщенной теоремы Адамяна--Арова--Крейна  нам понадобится  обобщение теоремы Нехари
для операторов Ганкеля на группах. Отметим, что теорема Нехари для ганкелевых форм на абелевы группы была перенесена
 Вонгом \cite{Wang}.

 Пусть $k$ --- функция на $X_+$. \textit{Ганкелевой формой} на $X_+$ c ядром $k$ называют комплексную
билинейную форму вида
$$
A(a,b)=\sum\limits_{\mu,\nu\in X_+}k(\mu\nu)a(\mu)b(\nu),\eqno(1.1)
$$
где $a,b$ ---  функции  на $X_+$.

{\bf  Теорема Нехари--Вонга \cite{Wang}}. {\it Ганкелева форма (1.1)
ограничена на $\ell^2(X_+)$ тогда и только тогда, когда ее ядро имеет вид
$k(\chi)=\widehat{\varphi}(\overline\chi),\ \chi\in X_+,$
 где $\varphi\in L^{\infty}(G)$.  При этом
$\|\varphi\|_{\infty}\leq M,$  где $M$ --- константа
ограниченности формы $A.$
В частности,  норма формы $A$ равняется
$\|\varphi\|_{\infty}$ для некоторой функции $\varphi$, удовлетворяющей
указанным выше условиям}.

 Доказательство. \footnote{Доказательство необходимости в \cite{Wang} содержит пробел, который мы устраняем.}
  Достаточность непосредственно следует из
вида ядра и неравенства $\|\varphi\|_{\infty}\leq M$.
Действительно, пусть функции $a$ и $b$ на $X_+$ имеют единичную $\ell^2(X_+)$-норму. Тогда по
теореме Планшереля функции $\alpha,\beta\in H^2(G)$, определяемые
с помощью сходящиеся рядов
$$\alpha(x)=\sum\limits_{\mu\in X_+}a(\mu)\mu(x),\ \beta(x)=\sum\limits_{\nu\in X_+}b(\nu)\nu(x),$$
имеют единичную $L^2$-норму. Поэтому
$$|A(a,b)|=\left|\sum\limits_{\mu,\nu\in X_+}\int\limits_G\varphi(x)(\mu\nu)(x)\,dxa(\mu)b(\nu)\right|=$$
$$=\left|\int\limits_G\varphi(x)\alpha(x)\beta(x)\,dx\right|\leq\|\varphi\|_{\infty}\|\alpha\|_2\|\beta\|_2\leq M.$$

Докажем теперь необходимость. Пусть $A$ --- ганкелева форма на $\ell^2(X_+)$ с ядром $k$,
ограниченная константой $M$. Покажем, что для любой  функции $g\in
H^1(G), \widehat g(1)\ne 0$  ряд
$$
\sum\limits_{\chi\in X_+}k(\chi)\widehat{g}(\chi)\eqno (1.2)
$$
сходится. Действительно, (см. \cite{Rud}, теорема 8.4.4) существуют такие
$\alpha,\beta\in H^2(G),$ что $g=\alpha\beta$ и
$\|\alpha\|_2=\|\beta\|_2=\|g\|_1^{1/2}.$ Обозначим
преобразования Фурье функций $\alpha$ и $\beta$ через $a$ и
$b$ соответственно. Тогда по теореме Планшереля
$$
\sum\limits_{\mu\in X_+}|a(\mu)|^2=\sum\limits_{\nu\in X_+}|b(\nu)|^2=\|g\|_1,
$$
$\widehat g=a\ast b,$
 а потому ряд (1.2) сводится к $A(a,b)$ и, следовательно, сходится.
В случае $\widehat g(1)\ne 0$ определим $\Lambda(g)$ как сумму
этого ряда. Имеем
$$
|\Lambda(g)|=|A(a,b)|\leq M\|a\|_2\|b\|_2= M\|g\|_1.\eqno (1.3)
$$

 Пусть теперь $g\in
H^1(G), \widehat g(1)=0.$ Тогда существует представление $g=g_1-g_2$, где  $g_i\in
H^1(G), \widehat g_i(1)\ne 0$ (например, $g=(g+1)-1$). В этом случае положим
$\Lambda(g)=\Lambda(g_1)-\Lambda(g_2).$ Легко видеть, что это определение корректно, и мы получаем линейный
функционал на $H^1(G)$. Более того,
из представления  $g=(g+\varepsilon)-\varepsilon$  ($\varepsilon$ --- любое положительное число)  с учетом (1.3)
 вытекает, что
$$
|\Lambda(g)|\leq |\Lambda(g+\varepsilon)|+|\Lambda(\varepsilon)|\leq M\|g\|_1+\varepsilon (M+|\Lambda(1)|),
$$
откуда в силу произвольности $\varepsilon$ следует, что $\Lambda$ является ограниченным функционалом  с нормой, не
превосходящей $M$. По теореме Хана-Банаха существует сохраняющее
норму расширение $\Lambda$ на $L^1(G)$, которое, как известно, задается некоторой функцией
$\varphi\in L^{\infty}(G)$ по формуле
$$
\Lambda(g)=\int\limits_G\varphi(x)g(x)\,dx, g\in L^1(G),\eqno (1.4)
$$
 причем $\|\varphi\|_{\infty}=\|\Lambda\|\leq M.$ Следовательно,  $\|\varphi\|_{\infty}\leq \|A\|.$
Из (1.4) при $g=\chi$ получаем равенство $k(\chi)=\widehat{\varphi}(\overline\chi).$
В свою очередь, из него  следует, что  $\|A\|\leq \|\varphi\|_{\infty}$ (см. выше доказательство достаточности).$\Box$

Версия теоремы Нехари--Вонга для операторов Ганкеля выглядит следующим образом.

\textbf{Теорема 1.5} (Нехари--Вонга).
 \textit{Ограниченный оператор $H:H^2(G)\to  H^2_-(G)$ имеет вид $H_\varphi$ для некоторой функции $\varphi\in L^\infty(G),$
  если и только если
$$
HS_\chi=P_-S_\chi H\  \forall \chi\in X_+,
$$
где $S_\chi f:=\chi f\ (f\in L^2(G))$. Кроме того, $\|H_{\varphi}\|={\rm dist}_{L^{\infty}}(\varphi,H^{\infty}(G)).$}

Для доказательства можно рассмотреть билинейную форму $(f,g)\mapsto\langle Hf,\overline g\rangle$ на $H^2(G)\times H^2(G)$
 и применить предыдущую теорему из  \cite{Wang} к соответствующей билинейной ганкелевой форме на
  $\ell^2(X_+)\times \ell^2(X_+)$ (преобразование Фурье изоморфно отображает $\ell^2(X_+)$ на $H^2(G)$),
  при этом получаем, что $\|H_{\varphi}\|=\|\varphi_0\|_\infty$ для некоторой функции $\varphi_0\in \varphi+H^\infty(G);$
   подробности см. в \cite{MD3}. \footnote{ В формулировке теоремы в \cite{MD3} имеется неточность;
    исправленный вариант появился в \cite{Indag}.}

Наконец, заметим, что $H_{\varphi}=H_{\varphi+\psi},$ если и только если $\psi\in
H^{\infty}(G).$ Поэтому
 для любых $\varphi\in L^\infty, \psi\in H^{\infty}(G)$ и $f\in H^2(G)$ имеем
$$
\|H_\varphi f\|_2=\|H_{\varphi+\psi}f\|_2=\|P_-((\varphi+\psi)f)\|_2\leq\|(\varphi+\psi)f\|_2\leq\|\varphi+\psi\|_{\infty}\|f\|_2.
$$
Следовательно,
$$
\|H_{\varphi}\|\leq\inf\limits_{\psi\in H^{\infty}(G)}\|\varphi+\psi\|_{\infty}={\rm dist}_{L^{\infty}}(\varphi,H^{\infty}(G)).
$$
А с другой стороны, так как $\|H_{\varphi}\|=\|\varphi_0\|_\infty$ для некоторой функции $\varphi_0\in \varphi+H^\infty(G),$
 то $\|H_{\varphi}\|\geq {\rm dist}_{L^{\infty}}(\varphi,H^{\infty}(G)),$
что и завершает
доказательство теоремы.$\Box$

\

\begin{center}
\S 2. Компактность и принадлежность классам Шаттена--фон Неймана операторов Ганкеля над компактными абелевыми группами
\end{center}

\

В этом разделе, в частности,  обобщаются и дополняются результаты работы \cite{MD4}.

Следующая теорема --- это  фактически  \cite[теорема 1.2]{YCG}.\footnote{  Определение ганкелева оператора в \cite{YCG}
отличается от данного выше, но доказательство легко приспособить к нашему случаю.}

\textbf{Теорема 2.1}   \cite{YCG}.
1)  \textit{На  $H^2(G)$ существует нетривиальный компактный оператор Ганкеля, если и только если  $X$ содержит наименьший
 положительный элемент} $\chi_1.$

2) \textit{Оператор $H_\varphi$ ($ \varphi\in L^\infty(G)$) компактен, если и только если  $P_-\varphi\in K(G)$},
\textit{где $K(G)=\mathrm{Cl}_{L^\infty}(\mathrm{span}\{\overline{\chi_1}^n: n\in \mathbb{N}\})$}.

Отметим, что первое утверждение теоремы 2.1 будет также установлено в ходе доказательства утверждения 1) теоремы 2.2,
 в которой дается  описание компактных операторов Ганкеля в $H^2(G).$  Эта часть теоремы 2.2 фактически содержится в
 доказательстве теоремы 1.2 из \cite{YCG}; для полноты  мы приводим подробное доказательство. Ниже мы используем
 обозначение $H^2({\chi_1})$ для подпространства в $H^2(G)$ с ортонормированным базисом $\{\chi_1^n: n\in\mathbb{ Z}_+\},$ $H^2(\chi_1)^\bot$ будет обозначать ортогональное дополнение
$H^2(G)\ominus H^2(\chi_1).$

\textbf{Теорема 2.2}. \textit{Пусть $X$ содержит наименьший положительный элемент $\chi_1,$ и пусть оператор
$H_{\varphi}$ ($\varphi\in L^\infty(G)$) компактен. Тогда}

1) \textit{$P_-\varphi$ разлагается в ряд Фурье по ортонормированной системе $\{\overline{\chi_1}^n:n\in \mathbb{N}\},$}
$$
P_-\varphi=\sum_{n=1}^\infty c_n \overline{\chi_1}^n.
$$

2) \textit{Сужение
$H_{\varphi}|H^2(\chi_1)^\bot$ равно нулю, а оператор $H_{\varphi}|H^2(\chi_1)$ унитарно эквивалентен классическому оператору
 Ганкеля $H_{\varphi_1}: H^2(\mathbb{T})\to  H^2_-(\mathbb{T})$ с символом
$$
\varphi_1(z)=\sum_{n=1}^\infty c_n \overline{z}^n,
$$
для которого $\overline{\varphi_1}\in VMOA(\mathbb{T}).$}

\textit{Обратно,  оператор $H_{\varphi},$ обладающий свойствами, описанными в \textrm{2)}, компактен.}

Доказательство. 1) Пусть  оператор $H_{\varphi}$ компактен и отличен от нуля. Обозначим через $\tau(x)$   оператор
 сдвига  на элемент $x\in G,$   $\tau(x)f(y)=f(xy),$ действующий в пространстве $L^2(G)$ и его подпространстве $H^2(G).$ Так как
  $\tau(x)\chi=\chi(x)\chi \ (\chi\in X),$ то $\tau(x)$ коммутирует с $P_\pm;$ в частности, $\tau(x)H_\varphi\tau(x^{-1})
  =H_{\tau(x^{-1})\varphi},$ а потому последний оператор компактен вместе с  $H_{\varphi}.$ Далее, так как
   $\varphi\notin H^\infty(G)$, то $\widehat\varphi(\chi)\neq 0$ для некоторого характера $\chi\in X_-.$ В силу
    непрерывной (и линейной) зависимости оператора  $H_{\varphi}$ от $\varphi\in L^\infty(G),$ \textit{для любого}
    такого характера компактный оператор
 $$
 \int\limits_G\chi(x)H_{\tau(x^{-1})\varphi}\,dx
 $$
(интеграл Бохнера) равен $H_{\chi\ast\varphi}=\widehat\varphi(\chi)H_\chi,$ откуда следует компактность оператора $H_\chi$
 (здесь мы воспользовались одной остроумной идеей из \cite{YCG}).
Далее, вычисляя $S_{\overline\chi} H_\chi\xi\ (\xi\in X_+)$ (черта обозначает комплексное сопряжение), получаем, что компактный оператор $S_{\overline\chi} H_\chi:
H^2(G)\to L^2(G)$ является проектором на подпространство, порожденное множеством характеров $[1,\overline{\chi})=[1,\chi^{-1}).$ Значит, это множество конечно, а потому содержит наименьший положительный
 элемент $\chi_1$. Кроме того, из конечности этого промежутка вытекает, что  $\chi\in X^i$ и $\chi=\overline{\chi_1}^n,
 \ n\in \mathbb{N}$  (см.   \cite[определение 3, теорема 2]{SbMath}). Следовательно, имеет место разложение в ряд Фурье
 $P_-\varphi=\sum_{n=1}^\infty c_n\overline{\chi_1}^n$.

2) В силу лемм 1.4 и 1.1, при $\xi\in X_+\setminus X^i$ имеем $H_{\varphi}\xi=P_-(P_-\varphi)\xi)=0,$ так как
$\overline{\chi_1}^n\xi\in X_+\ (n\in \mathbb{Z}_+)$  по лемме 1.1. Значит, $H_{\varphi}|H^2(\chi_1)^\bot=0.$

Далее через $U_+ (U_-)$ будут обозначаться унитарные операторы   из $H^2(\mathbb{T})$ в  $H^2(\chi_1)$
(соответственно, из $H^2_-(\mathbb{T})$  в $H^2(\overline{\chi_1})$), определяемые равенствами
 $U_+z^m=\chi_1^m$ при $m\in \mathbb{Z}_+$  и  $U_-z^m=\chi_1^m$ при $m\in \mathbb{Z}, m<0$ (т.~е. $U_\pm f=f\circ\chi_1$).
 Справедливо равенство $H_{\varphi}|H^2(\chi_1)=U_-^{-1}H_{\varphi_1}U_+,$   т.~е. коммутативна следующая диаграмма:
\[
\begin{CD}
H^2(\mathbb{T})@>U_+>>H^2(\chi_1)\\
@VV{H_{\varphi_1}}V               @VV{H_{\varphi}|H^2(\chi_1)}V\\
H^2_-(\mathbb{T})@>U_->>H^2(\overline{\chi_1}).\\
\end{CD}
\]
 В самом  деле, разложение $P_-\varphi$ в ряд по степеням $\overline{\chi_1},$  установленное в 1), позволяет
 легко проверить, что обе части этого равенства совпадают на элементах базиса $\{\chi_1^n: n\in \mathbb{Z}_+\}.$
 Таким образом, оператор $H_{\varphi_1}$ также компактен. Кроме того, $\varphi_1\in L^2(\mathbb{T})$ так как $(c_n)
 \in \ell^2.$ Поэтому ($P_-^1$ обозначает ортопроектор из  $L^2(\mathbb{T})$ на  $H^2_-(\mathbb{T})$) $\varphi_1=P_-^1
 \varphi_1\in VMO(\mathbb{T})$  (см., например, \cite[теорема 5.8]{Pel}). Поскольку функция $\overline{\varphi_1}$
 аналитична, отсюда следует, что $\overline{\varphi_1}\in VMOA(\mathbb{T}).$

 Справедливость обратного утверждения легко следует из уже упоминавшегося описания компактных ганкелевых операторов
  на группе окружности. $\Box$

Теорема 2.2 позволяет получить обобщения классических теорем   Хартмана и Пеллера.

Как отмечалось выше, обобщение теоремы Хартмана на случай компактных абелевых групп было дано в \cite{YCG}. Следующая теорема
  дает иной,
 чем в  \cite{YCG}, критерий компактности операторов Ганкеля.

\textbf{Теорема 2.3}. \textit{Пусть $X$ содержит наименьший положительный элемент $\chi_1,$ и $\varphi\in L^\infty(G).$
 Тогда оператор $H_{\varphi}$ компактен, если и только если
$P_-\varphi=\varphi_1\circ \chi_1,$ где $\varphi_1\in \mathrm{VMO}(\mathbb{T}).$}

Доказательство.   Из теоремы 2.2 сразу следует, что оператор $H_{\varphi}$ компактен тогда и только тогда, когда
компактен оператор $H_{\varphi_1}.$   Это, в свою очередь, равносильно условию $\varphi_1\in \mathrm{VMO}(\mathbb{T})$
  (см., например, \cite[c. 52, теорема 5.8]{Pel}).
Осталось воспользоваться равенством  $P_-\varphi=\varphi_1\circ \chi_1.$$\Box$

  Приведем теорему, которая переносит на группы важную теорему В.~ В.~ Пеллера \cite{Pel}.  Определение пространств
   Бесова можно также найти в  \cite{Pel}.

\textbf{Теорема 2.4.} \textit{Пусть $X$ содержит наименьший положительный элемент $\chi_1,$ и $\varphi\in L^\infty(G),$ $p>0.$
 Тогда оператор $H_{\varphi}$ принадлежит классу Шаттена--фон Неймана $\mathbf{S}_p(H^2(G),H^2_-(G)),$ если и только если
$P_-\varphi=\varphi_1\circ \chi_1,$ где $\varphi_1$ принадлежит пространству Бесова  $B_p^{1/p}.$}

Доказательство.  Снова воспользовавшись теоремой 2.2 и ее доказательством, получаем, что  оператор $H_{\varphi}$
принадлежит  $\mathbf{S}_p(H^2(G),H^2_-(G)),$ если и только если оператор $H_{\varphi}|H^2(\chi_1)$
принадлежит  $\mathbf{S}_p(H^2(\chi_1),H^2(\overline{\chi_1}))$ (в ходе доказательства теоремы 2.2 было показано,
 что $H_{\varphi}|H^2(\chi_1)$ отображает $H^2(\chi_1)$ в  $H^2(\overline{\chi_1})$). В свою очередь, по этой же теореме
последнее равносильно тому, что
$H_{\varphi_1}$   принадлежит  $\mathbf{S}_p(H^2(\mathbb{T}),H^2_-(\mathbb{T})).$ Так как $P_-^1\varphi_1 =\varphi_1,$
 то,  согласно упоминавшимся результатам В.~В.~ Пеллера  (\cite[с. 314, следствие 1.2; c. 323, следствие 2.2; с. 328,
  следствие 3.2]{Pel}), это равносильно условию $\varphi_1\in B_p^{1/p}$  (указанные следствия применимы, так как в силу
  теоремы 2.2 $\varphi_1\in\overline{ \mathrm{VMOA}(\mathbb{T})}\subset \mathrm{BMO}(\mathbb{T})).$$\Box$

\textbf{Следствие 2.5} \cite{MD4}. \textit{Пусть   $\varphi\in L^\infty(G)$. Оператор  $H_\varphi$  является ядерным,
 если и только если найдутся последовательности
комплексных чисел $(c_n)$ и $(\lambda_n),$ такие, что $|\lambda_n|<1, \sum_{n=1}^\infty|c_n|/(1-|\lambda_n|)<\infty$  и
\[
P_-\varphi= \sum\limits_{n=1}^\infty\frac{c_n}{1-\lambda_n\overline{\chi_1}}.
 \]
 }

Доказательство. В силу теоремы 2.4 оператор $H_\varphi$  является ядерным, если и только если
$P_-\varphi=\varphi_1\circ \chi_1,$ где $\varphi_1\in B_1^{1},$ или, что равносильно, $\overline{\varphi_1}\in \overline{B_1^{1}}
=B_1^{1}.$ Поскольку $\overline{\varphi_1}$ аналитична, последнее в точности означает, что  $\overline{\varphi_1}$
 принадлежит подпространству $(B_1^{1})_+$ пространства $B_1^{1},$ состоящему из аналитических функций. Но известно
 \cite[c. 319, теорема 1.4 ]{Pel}, что это, в свою очередь, равносильно тому, что
найдутся последовательности
комплексных чисел $(c'_n)$ и $(\lambda'_n),$ такие, что $|\lambda'_n|<1, \sum_{n=1}^\infty|c'_n|/(1-|\lambda'_n|)<\infty$  и
\[
\overline{\varphi_1(\zeta)}= \sum\limits_{n=1}^\infty\frac{c'_n}{1-\lambda'_n\zeta},
 \]
что и завершает доказательство следствия.$\Box$

\

\begin{center}
\S 3. Конечные произведения Бляшке и конечномерность операторов Ганкеля над компактными абелевыми группами
\end{center}

\

\textbf{Теорема 3.1}. \textit{Пусть $X$ содержит наименьший положительный элемент $\chi_1,$ и $\varphi\in L^\infty(G)$.
 Тогда оператор $H_{\varphi}$ имеет конечный ранг, если и только если
$P_-\varphi=R\circ \chi_1,$ где $R$  есть рациональная функция, все полюсы которой лежат в $\mathbb{D},$ и при этом}
$$
\mathrm{rank }H_{\varphi}=\mathrm{deg} R.
$$

Доказательство. В силу теоремы 2.2 $\mathrm{rank}H_{\varphi}=\mathrm{rank}H_{\varphi_1}.$ Значит, $\mathrm{rank}H_{\varphi}$
 конечен тогда и только тогда, когда конечен $\mathrm{rank}H_{\varphi_1}.$  В силу классической теоремы Кронекера последнее
 условие выполняется тогда и только тогда, когда $P_-^1{\varphi_1}=R,$ где $R$  есть рациональная функция, все полюсы которой
  лежат в $\mathbb{D},$ и при этом $\mathrm{rank }H_{\varphi_1}=\mathrm{deg} R$  (как и выше, $P_-^1$ обозначает ортопроектор
  из  $L^2(\mathbb{T})$ на  $H^2_-(\mathbb{T})$). Поэтому $P_-^1{\varphi_1}=\varphi_1$ в силу  определения функции $\varphi_1.$
   Значит,  $\mathrm{rank }H_{\varphi}=\mathrm{deg} R$ и из определения функции $\varphi_1$ сразу следует, что
    $P_-\varphi=\varphi_1\circ \chi_1=R\circ \chi_1.$ Обратно, если $P_-\varphi=R\circ \chi_1,$ то
$\varphi_1\circ \chi_1=R\circ \chi_1.$ Заметим, что множество значений характера  $\chi_1$ есть
единичная окружность $\mathbb{T},$ так как это множество есть нетривиальная связная подгруппа группы  $\mathbb{T}.$
 Значит, $\varphi_1=R,$ и обратное утверждение теоремы также следует из теоремы Кронекера.$\Box$

\textbf{Определение 3.2.} Следуя \cite{MD4}, внутреннюю функцию  $\Theta$ (т.~е. такую функцию $\Theta\in H^2(G),$  что $|\Theta|=1$ п.в., см.
 \cite{Rud}) назовем \textit{конечным произведением Бляшке на группе  $G,$} если пространство
$$
K_\Theta:=H^2(G)\ominus\Theta H^2(G)
$$
 конечномерно. При этом мы положим ${\rm deg} \Theta:={\rm dim} K_\Theta$.

 Для описания конечных произведений Бляшке на группе  $G$ нам понадобятся обобщения двух классических лемм из \cite{MD4}.
 Для удобства читателя мы приводим их доказательство. Всюду ниже   $S_\chi:H^2(G)\to H^2(G), f\mapsto \chi f$ --- \textit{оператор умножения на характер }  $\chi\in X_+$.

 \textbf{Лемма 3.3} (обобщенная лемма Бёрлинга).  1) \textit{Любое число $\lambda\in\mathbb{D}$  является собственным значением оператора $S_\chi^*\ (\chi\in X_+,\ \chi\ne 1),$ и  отвечающие ему собственные функции имеют в точности вид $g/(1-\lambda\chi)$,  где $g\in H^2(G),\ g\ne 0$, причем  преобразование Фурье $\widehat g$  сосредоточено на промежутке $[1,\chi)$.}

2) \textit{Пусть в группе $X$ имеется наименьший положительный элемент $\chi_1$. Тогда точечный спектр  оператора $S_{\chi_1}^*$ равен $\mathbb{D}$, и при всех $\lambda\in \mathbb{D}$ и натуральных  $n$}
$$
{\rm Ker}(S_{\chi_1}^*-\lambda I)^n={\rm span}\left\{\frac{\chi_1^{k-1}}{(1-\lambda\chi_1)^k}: 1\leq k\leq n\right\}.
$$

Доказательство. 1)   Заметим, что  $\sigma(S_\chi^*)\subseteq\overline{\mathbb{D}}$, так как $\|S_\chi^*\|=1$. Если число $\lambda\in\overline{\mathbb{D}}$ есть собственное значение оператора $S_\chi^*\ (\chi\in X_+,\ \chi\ne 1)$ с собственной функцией $f,$ то, учитывая, что $S_\chi^*=P_+S_{\overline \chi}$, и  приравнивая коэффициенты Фурье функций $S_\chi^* f$ и  $\lambda f (f\in H^2(G),\ f\ne 0)$, получаем, что $\widehat f(\eta\chi)=\lambda \widehat f(\eta)\ (\eta\in X_+)$, что равносильно тому, что преобразование Фурье функции $f(1-\lambda\chi)\chi^{-1}=f\chi^{-1}-\lambda f$ сосредоточено на  $X_-$. Если мы положим $g:=f(1-\lambda\chi)$, то это значит, что   $\widehat g$ сосредоточено на $[1,\chi)$. При $\lambda\in\mathbb{D}$ получаем требуемый вид собственных функций, отвечающих $\lambda$.

2) При $\chi=\chi_1$ из только что доказанного утверждения 1) леммы следует, что функция $\widehat g$ сосредоточена на $\{1\}$, т.~е. $g$ есть константа, отличная от $0$. Если мы предположим, что $\lambda\in\mathbb{T}$ есть собственное значение оператора  $S_{\chi_1}^*$, то по доказанному выше  $g=f(x)(1-\lambda\chi_1(x))$ при всех $x\in G,$ что невозможно, так как $\chi_1$ принимает значение $1/\lambda$\ ($\chi_1(G)=\mathbb{T}$ как нетривиальная связная подгруппа группы  $\mathbb{T}$).

Докажем последнее утверждение леммы индукцией по $n$. При  $n=1$  оно уже доказано. Предположим, что оно верно при некотором $n$, и рассмотрим функцию $f\in {\rm Ker}(S_{\chi_1}^*-\lambda I)^{n+1}$. Так как $(S_{\chi_1}^*-\lambda I)f\in{\rm Ker}(S_{\chi_1}^*-\lambda I)^n$, то по предположению
$$
(S_{\chi_1}^*-\lambda I)f=\sum\limits_{k=1}^nc_k\frac{\chi_1^{k-1}}{(1-\lambda\chi_1)^k}.
$$
Но легко проверить, что
$$
\frac{\chi_1^{k-1}}{(1-\lambda\chi_1)^k}=(S_{\chi_1}^*-\lambda I)\frac{\chi_1^{k}}{(1-\lambda\chi_1)^{k+1}}.
$$
Поэтому
$$
(S_{\chi_1}^*-\lambda I)\left(f-\sum\limits_{k=1}^n c_k\frac{\chi_1^{k}}{(1-\lambda\chi_1)^{k+1}}\right)=0,
$$
откуда по доказанному выше найдется такая константа  $c_0$, что
$$
f-\sum\limits_{k=1}^n c_k\frac{\chi_1^{k}}{(1-\lambda\chi_1)^{k+1}}=\frac{c_0}{1-\lambda\chi_1},
$$
что и требовалось доказать.$\Box$

 \textbf{Лемма 3.4}  (ср. \cite[с. 216, лемма 2.4.4]{Nik1}).  \textit{Пусть в группе $X$ имеется наименьший положительный элемент $\chi_1$. Если  $\Theta$ --- конечное произведение Бляшке на группе $G,$  то для некоторого конечного множества $E\subset \mathbb{D}$ и некоторой функции $k:E\to \mathbb{N}$ имеем}
$$
K_\Theta={\rm span}\left\{\frac{\chi_1^{k}}{(1-\lambda\chi_1)^{k+1}}: \lambda\in E, 1\leq k\leq k(\lambda) \right\}.
$$

Доказательство. Поскольку оператор  $S_{\chi_1}^*|K_\Theta$ действует в пространстве $K_\Theta,$ это следует из найденного в лемме 3.3  вида корневых подпространств  оператора $S_{\chi_1}^*$ и того факта, что конечномерное пространство  есть прямая сумма корневых подпространств действующего в нем оператора.$\Box$

Нам понадобится также следующая элементарная

\textbf{Лемма 3.5.} \textit{Пусть  $L, L_1$ --- замкнутые подпространства гильбертова пространства и  $L_1\subset L.$ Тогда
}$$
(L_1\oplus L^\bot)^\bot=L\ominus L_1.
$$

 \textbf{Теорема 3.6.} \textit{На группе  $G$ существует конечное произведение Бляшке тогда и только тогда, когда  $X$ содержит наименьший положительный элемент $\chi_1.$ При этом конечные произведения Бляшке на $G$ имеют в точности вид  $\Theta=B\circ \chi_1,$ где $B$ --- конечное произведение Бляшке на группе $\mathbb{T},$ и ${\rm deg} \Theta={\rm deg} B$.}

Доказательство. Если $\Theta$ есть конечное произведение Бляшке на  $G,$ то образ ганкелева оператора $H_{\overline{\Theta}}$ содержится в $K_\Theta$ ($\Theta H^2(G)\subset H^2(G)$  для любой внутренней функции $\Theta,$ а  потому $\mathrm{Im} (H_{\overline{\Theta}})\subset H^2_-(G)\subset (\Theta H^2(G))^\bot$). Следовательно, оператор $H_{\overline{\Theta}}$ компактен и нетривиален (в случае его тривиальности  $\overline{\Theta}H^2(G)\subset H^2_-(G),$ что невозможно, так как $1\in \overline{\Theta}H^2(G)$), и по теореме  2.1 $X$ содержит наименьший положительный элемент.

Далее, пусть $\Theta:=B\circ \chi_1,$ где $B$ --- конечное произведение Бляшке на  $\mathbb{T}.$
Так как $X_+\setminus X^i$ есть ортонормированный базис пространства  $H^2(\chi_1)^\bot,$ а $\Theta=B\circ \chi_1=U_+(B)\in H^2(\chi_1),$  то по утверждению 2) леммы 1.1 $\Theta (H^2(\chi_1))^\bot=H^2(\chi_1)^\bot.$  Кроме того, $\Theta H^2(\chi_1)\subset H^2(\chi_1).$
Следовательно, применяя лемму 3.5 и теорему 2.2, получаем
$$
K_\Theta= (\Theta H^2(G))^\bot=
$$
$$
(\Theta H^2(\chi_1)\oplus \Theta (H^2(\chi_1))^\bot)^\bot=(\Theta H^2(\chi_1)\oplus H^2(\chi_1)^\bot)^\bot
$$
$$
 =H^2(\chi_1)\ominus \Theta  H^2(\chi_1)=U_+(H^2(\mathbb{T})\ominus BH^2(\mathbb{T}))=U_+(K_B),
$$
 а потому $\Theta$ --- конечное произведение Бляшке  и ${\rm deg} \Theta={\rm deg} B$.

Наконец, в силу леммы 3.4 и описания классических пространств $K_B$ \cite[с. 216, лемма 2.4.4]{Nik1}, для любого конечного произведения Бляшке $\Theta$ на $G$  найдется такое конечное произведение Бляшке $B_1$ на группе $\mathbb{T},$ что
$$
K_\Theta=K_{B_1}\circ\chi_1:=\{f\circ\chi_1:f\in K_{B_1}\}=K_{B_1\circ\chi_1}.
$$
  Следовательно, для внутренних функций $\Theta$ и $\Theta_1:=B_1\circ\chi_1$ имеем  $\Theta f=\Theta_1$ и $\Theta_1 g=\Theta$ для некоторых $f,g\in H^2(G)$. Из этих равенств следует, что $|f|=|g|=1$ и $fg=1.$ Таким образом,
$g=\overline{f}\in H^2(G)\cap \overline{H^2(G)}=\{\mathrm{const}\}.$ Значит, $g(x)=c\in \mathbb{T}$ и
$\Theta =c\Theta_1=B\circ\chi_1,$ где $B:=cB_1.$$\Box$

 \textbf{Теорема 3.7.} \textit{Пусть $\varphi\in H^{\infty}(G).$ Условие
  $\Theta P_-\varphi\in H^{\infty}(G)$ для некоторого конечного произведения Бляшке $\Theta$  необходимо и достаточно для
   конечномерности оператора $H_{\varphi}.$
}

Доказательство. Пусть оператор $H_{\varphi}=H_{P_-\varphi}$ имеет конечный ранг.
 Если $P_-\varphi=\varphi_1\circ \chi_1$  как в теореме 2.2, то по этой теореме оператор $H_{\varphi_1}$
также конечномерен. В силу классического результата (см. \cite[глава 1, следствие 3.3]{Pel}) найдется такое конечное произведение Бляшке $B$ на группе  $\mathbb{T},$ что $B\varphi_1\in H^\infty(\mathbb{T}).$ Следовательно, $(B\circ \chi_1)(\varphi_1\circ \chi_1)\in H^\infty(G)$ и для завершения  доказательства необходимости осталось воспользоваться теоремой 3.6.

Обратно, если $\Theta P_-\varphi\in H^{\infty}(G)$ для некоторого конечного произведения Бляшке $\Theta$, то  $H_{\Theta P_-\varphi}=0$, т.~е. $H_{\varphi}(\Theta f)=0$ при всех $f\in H^2(G)$. Значит, $\Theta H^2(G)\subseteq H_0,$ где $H_0:={\rm Ker}H_{\varphi},$ а, стало быть, $H_0^\bot\subseteq K_\Theta$. Таким образом, ${\rm dim}H_0^\bot<\infty$, а так как $H_{\varphi}(H^2(G))=H_{\varphi}(H_0^\bot)$, причем оператор $H_{\varphi}|H_0^\bot$ инъективен, то ${\rm rank}H_{\varphi}={\rm dim}H_0^\bot$.$\Box$

\textbf{Следствие 3.8.} \textit{Если  $\Theta$ --- конечное произведение Бляшке,  $h\in H^\infty(G)$ и $P_-\varphi=\overline \Theta h,$ то  ${\rm rank}H_{\varphi}<\infty$.}

\

\begin{center}
\S 4. Теорема Адамяна--Арова--Крейна для операторов Ганкеля над компактными абелевыми группами
\end{center}

\

Следующий результат обобщает классическую теорему  Адамяна--Арова--Крейна (см., например, \cite{Pel}). Как и выше,
доказательство использует редукцию к классическому случаю, описанную в теореме 2.2.

\textbf{Теорема 4.1}. \textit{Пусть $X$ содержит наименьший положительный элемент $\chi_1,$ $\varphi\in L^\infty(G)$ и
  оператор $H_{\varphi}$ компактен. Тогда для сингулярных чисел оператора $H_{\varphi}$ справедливы равенства ($n\geq 1$)
$$
s_n(H_{\varphi})=\min\{\|H_{\varphi}-H_n\|: \mbox{ оператор }  H_n \mbox{ ганкелев над } G, \mathrm{rank} H_n\leq n\}.
$$
}

Доказательство. В силу упоминавшейся теоремы Адамяна--Арова--Крейна,
$$
s_n(H_{\varphi_1})=\min\{\|H_{\varphi_1}-\Gamma_n\|: \mbox{ оператор }  \Gamma_n \mbox{ ганкелев над }
\mathbb{T}, \mathrm{rank} \Gamma_n\leq n\},\eqno(4.1)
$$
причем, в соответствии с  теоремой 2.2,
$$
s_n(H_{\varphi_1})=s_n(H_{\varphi}|H^2(\chi_1)).\eqno(4.2)
$$
С другой стороны, разлагая оператор $H_{\varphi}|H^2(\chi_1)$  в ряд Шмидта, получаем для $x_1\in H^2(\chi_1)$
$$
H_{\varphi}x_1=(H_{\varphi}|H^2(\chi_1))x_1=\sum_n s_n(H_{\varphi}|H^2(\chi_1))\langle x_1,e_n\rangle e'_n,
$$
где $(e_n),$ $(e'_n)$--- ортонормированные базисы пространств $H^2(\chi_1)$  и  $H^2(\overline{\chi_1})\subset H^2_-(G)$
соответственно (в ходе доказательства теоремы 2.2 было показано, что $H_{\varphi}: H^2(\chi_1)\to H^2(\overline{\chi_1})$).

Далее, для любого $x\in H^2(G)$ имеем $x=x_1+x_2,$ где $x_1\in H^2(\chi_1),$  $x_2\in H^2(\chi_1)^\bot:=H^2(G)
\ominus H^2(\chi_1),$ причем  $H_{\varphi}x_2=0$ по теореме 2.2. Поэтому
$$
H_{\varphi}x=H_{\varphi}x_1=\sum_n s_n(H_{\varphi}|H^2(\chi_1))\langle x,e_n\rangle e'_n
$$
 ($\langle x_2,e_n\rangle=0$). Отсюда следует, что
$$
s_n(H_{\varphi}|H^2(\chi_1))=s_n(H_{\varphi}).\eqno(4.3)
$$
Объединяя равенства (4.3), (4.2) и (4.1), получаем
$$
s_n(H_{\varphi})=\min\{\|H_{\varphi_1}-\Gamma_n\|:   \Gamma_n \mbox{ ганкелев на группе } \mathbb{T}, \mathrm{rank}
 \Gamma_n\leq n\}.\eqno(4.4)
$$
В силу теоремы 2.2 $H_{\varphi_1}=U_-^{-1}(H_{\varphi}|H^2(\chi_1))U_+.$ С учетом унитарности операторов $U_+, U_-,$ имеем
$$
\|H_{\varphi_1}-\Gamma_n\|=\|H_{\varphi}|H^2(\chi_1)-U_-\Gamma_nU_+^{-1}\|=\|H_{\varphi}|H^2(\chi_1)-\widetilde{H_n}\|,
$$
где $\widetilde{H_n}:=U_-\Gamma_nU_+^{-1}:H^2(\chi_1)\to H^2(\overline{\chi_1}),$  и $\mathrm{rank}\widetilde{H_n}=\mathrm{rank}
 \Gamma_n\leq n.$

Определим оператор $H_n:H^2(G)\to H^2_-(G)$ с помощью равенств
$$
H_n|H^2(\chi_1):=\widetilde{H_n},\quad H_n|H^2(\chi_1)^\bot:=0.
$$
Ясно, что $\mathrm{rank}H_n=\mathrm{rank}\widetilde{H_n}\leq n,$  и
$$
\|H_{\varphi_1}-\Gamma_n\|=\|H_{\varphi}-H_n\|.\eqno(4.5)
$$
Покажем, что оператор $H_n$ ---  ганкелев.
С этой целью  проверим, что  $H_n$ удовлетворяет коммутационным соотношениям
$$
H_nS_\chi f=P_-S_\chi H_n f\quad (\chi\in X_+, f\in H^2(G))\eqno(4.6)
$$
(см. теорему 1.5), где $S\chi:L^2(G)\to L^2(G), h\mapsto \chi h$  --- оператор сдвига на характер.
Для этого достаточно рассмотреть следующие случаи.

1) $f\in H^2(\chi_1)^\bot.$ Так как  при $\chi\in X_+$ по лемме 1.1 $\chi(X_+\setminus X^i)\subset X_+\setminus X^i,$ то
$S_\chi$ отображает пространство $H^2(\chi_1)^\bot$ в себя (ортонормированный базис этого пространства, состоящий из
 характеров, содержится в $X_+\setminus X^i$). Поэтому в этом случае обе части в (4.6) равны нулю.

2)  $f\in H^2(\chi_1).$ В этом случае $H_nf=\widetilde{H_n}f.$ Поскольку $f$ разлагается в ряд по элементам $X_+\cap X^i,$
 можно считать, что $f=\chi_1^k, k\geq  0.$ Тогда (4.6) приобретает вид
$$
H_n(\chi\chi_1^k) =P_-(\chi \widetilde{H_n}\chi_1^k) \quad (\chi\in X_+, k\geq 0).\eqno(4.7)
$$
Возможны два случая.

а) $\chi\in X_+\cap X^i,$  т.~е. $\chi=\chi_1^j, j\geq  0.$ Тогда (4.7) равносильно равенству
$$
H_n(\chi_1^{j+k}) =P_-(\chi_1^j \widetilde{H_n}\chi_1^k) \quad (j, k\geq 0).
$$
Последнее равенство вытекает из коммутационных соотношений для оператора Ганкеля $\Gamma_n.$
В самом деле, пусть $L^2(\chi_1)$  есть  гильбертово пространство с ортонормированным базисом $\{\chi_1^m: m\in\mathbb{ Z}\},$
 и пусть $U:L^2(\mathbb{T})\to L^2(\chi_1),$
$z^m\mapsto\chi_1^m$ ($m\in\mathbb{ Z}$) есть  унитарный оператор, продолжающий $U_{\pm}.$
 Так как $\widetilde{H_n}\chi_1^m=U_-\Gamma_nU_+^{-1}\chi_1^m=U_-\Gamma_nz^m,$ то (4.7) принимает вид
$U_-\Gamma_nz^{k+j}=P_-(U(z^j)(U_-\Gamma_nz^k)).$ Далее, разлагая функцию $f\in L^2(\mathbb{T})$ в ряд Фурье, легко проверить,
что $U(z^j f)=U(z^j)U(f).$ Легко проверить также, что $P_-U=UP_-^1$
 (как и выше, $P_-^1$ обозначает ортопроектор из  $L^2(\mathbb{T})$ на  $H^2_-(\mathbb{T})$). Поэтому
$$
U_-\Gamma_nz^{k+j}=P_-U(z^j)(U_-\Gamma_nz^k)=P_-U(z^j\Gamma_nz^k)=UP_-^1(z^j\Gamma_nz^k).
$$
Данное же равенство вытекает из равенства $\Gamma_nz^{k+j}=P_-^1(z^j\Gamma_nz^k),$ являющегося следствием коммутационных
соотношений для $\Gamma_n.$

б)  $\chi\in X_+\setminus X^i.$  Тогда  по лемме 1.1 $\chi\chi_1^k\in X_+\setminus X^i\subseteq H^2(\chi_1)^\bot$  ($k\geq 0$).
Значит, в этом случае левая часть в (3.7) равна нулю. Далее, у нас $\widetilde{H_n}\chi_1^k\in H^2(\overline{\chi_1}).$
 Поэтому $\chi\widetilde{H_n}\chi_1^k\in H^2(G)$  в силу леммы 1.1, а потому и правая часть в (4.7) равна нулю.
 Итак, оператор $H_n$ ---  ганкелев.

 Наконец, по классической теореме Адамяна--Арова--Крейна существует ганкелев оператор $\Gamma_n^0$ на группе $\mathbb{T}$ такой,
  что $s_n(H_{\varphi_1})=\|H_{\varphi_1}-\Gamma_n^0\|.$ Пусть $H_n^0$ --- соответствующий $\Gamma_n^0$ ганкелев оператор
  на группе $G,$ конструкция которого была изложена выше. Тогда с учетом (4.5) имеем
 $$
 s_n(H_{\varphi})=s_n(H_{\varphi_1})=\|H_{\varphi_1}-\Gamma_n^0\|=\|H_{\varphi}-H_n^0\|\geq \inf\|H_{\varphi}-H_n\|\geq
  s_n(H_{\varphi})
 $$
(инфимум берется по всем ганкелевым операторам на группе $G$ ранга не выше $n$). $\Box$

Воспользуемся обозначением
$$
\mathcal{R}_n=\{R_n: R_n \mbox{ --- рациональная функция с полюсами в  } \mathbb{D}, \mathrm{deg}R_n\leq n \}.
$$

\textbf{Следствие 4.2.} \textit{Пусть $X$ содержит наименьший положительный элемент $\chi_1,$ $\varphi\in L^\infty(G)$ и
 оператор $H_{\varphi}$ компактен. Тогда}
$$
s_n(H_{\varphi})=\mathrm{dist}_{L^\infty}(\varphi, \mathcal{R}_n\circ\chi_1+H^\infty(G)),
$$
\textit{где  $\mathcal{R}_n\circ\chi_1:=\{R_n\circ\chi_1:R_n \in \mathcal{R}_n\},$ $n\in \mathbb{N}.$}

Доказательство. В силу теоремы 3.1,    ганкелев оператор $H_n$ на группе $G,$ удовлетворяющий условию $\mathrm{rank} H_n\leq n,$
 имеет вид  $H_n=H_{R_n\circ\chi_1},$ где $R_n\in \mathcal{R}_n.$ Следовательно, используя теорему 4.1 и теорему Нехари-Вонга,
  имеем
$$
s_n(H_{\varphi})=\min\{\|H_\varphi-H_{R_n\circ\chi_1}\|:R_n\in   \mathcal{R}_n\}=
$$
$$
\min\{\mathrm{dist}_{L^\infty}(\varphi-R_n\circ\chi_1,H^\infty(G)):R_n\in   \mathcal{R}_n\}=\mathrm{dist}_{L^\infty}
(\varphi, \mathcal{R}_n\circ\chi_1+H^\infty(G)).\Box
$$

\

\begin{center}
\S 5. Инвариантные подпространства $H^2(G)$
\end{center}

\

 Напомним (см., например, \cite{Rud}), что  подпространство $E$ пространства $H^2(G)$ называется \textit{инвариантным}, если $\chi E\subset E$ при всех $\chi\in X_+.$

\textbf{Определение 5.1.} Инвариантное подпространство $E$ пространства $H^2(G)$ назовем \textit{неприводящим}, если оно удовлетворяет условию $\chi E\subsetneqq  E$ при всех $\chi\in X_+\setminus\{1\}.$

Подпространство $E$ пространства $H^2(G)$ вида $\theta H^2(G),$ где $\theta$ --- внутренняя функция, является замкнутым, инвариантным и неприводящим,  так как в противном случае  при некотором $\chi\in X_+\setminus\{1\}$ выполнялось бы неверное равенство $\chi H^2(G)=H^2(G).$

Следующая теорема дополняет известное обобщение Хелсона и Лауденслегера теоремы Бёрлинга об инвариантных подпространствах \cite[теорема 8.5.3]{Rud}.

 \textbf{Теорема 5.2.} \textit{Каждое замкнутое  нетривиальное  неприводящее инвариантное подпространство $E$  пространства $H^2(G)$   имеет вид $\theta H^2(G),$ где $\theta$ --- внутренняя функция, тогда и только тогда, когда  $X$ содержит наименьший положительный элемент. При этом функция $\theta$ единственна с точностью до множителя из $\mathbb{T}.$
}

Доказательство. В \cite[8.5.4]{Rud} замечено, что если $X$ не содержит наименьшего положительного элемента, то  $E_0:=\{f\in H^2(G): \widehat{f}(1)=0\}$ есть  замкнутое  инвариантное подпространство пространства $H^2(G)$, которое не  имеет вида $\theta H^2(G),$ где $\theta$ --- внутренняя функция. В самом деле, пусть $E_0=\theta H^2(G)$ и $X$ не содержит наименьшего положительного элемента. Разложение в ряд Фурье функции $\theta\in E_0,$  имеет вид $\theta=\sum_{\chi\in  X_\theta} a_\chi \chi,$ где $X_\theta\subset  \{\chi\in X_+: \chi>1\}.$ Выберем характер $\xi\in X_+\setminus\{1\}$ так, что $\xi<\chi$ для всех $\chi\in X_\theta.$ Тогда $\xi\in E_0,$ но $\xi\notin \theta H^2(G),$ и мы получили противоречие. При этом $E_0$ неприводящее ($\chi E_0\ne  E_0$ при всех $\chi\in X_+\setminus\{1\}$, так как  $E_0$ содержит интервал $(1,\chi)\subset X_+$).

Пусть теперь $X$ содержит наименьший положительный элемент $\chi_1.$  Предположим, что замкнутое инвариантное подпространство $E$ является неприводящим.  Заметим, что $\chi E\subset \chi_1 E$ для любого $\chi\in X_+\setminus\{1\},$ поскольку в этом случае $\chi_1$ делит $\chi.$ Поэтому мы можем обобщить на наш случай известное доказательство теоремы Бёрлинга  (см., например, \cite[c. 9--10]{Nik1}).  Рассмотрим функцию $\theta\in E\cap (\chi_1 E)^\bot$ с условием нормировки $\|\theta\|_2=1.$ Тогда $\theta\bot\chi\theta$ для любого характера $\chi\in X_+\setminus\{1\},$ а потому
$$
\int\limits_G|\theta|^2\chi \,dx=0\quad (\chi\in X_+\setminus\{1\}).
$$
Взяв комплексное сопряжение, получим, что преобразование Фурье функции $|\theta|^2$ сосредоточено на $\{1\},$ т.е. $|\theta|^2=\mathrm{const},$ что вместе с условием нормировки  влечет равенство $|\theta|=1$ п.в. Так как $\theta\in E,$  то $\theta \mathrm{span} X_+\subset E,$ где $\mathrm{span} X_+$   есть множество всех аналитических полиномов на $G.$
Поскольку $\mathrm{span} X_+$  плотно в $H^2(G)$ по определению, получаем, что $\theta H^2(G)\subset E.$
Для доказательства обратного включения возьмем $f\in E\cap (\theta H^2(G))^\bot.$ Тогда $ f\bot \theta\chi,$ т.е. $\int_Gf\overline{\theta\chi}\, dx=0$ для любого характера $\chi\in X_+.$ С другой стороны,  для любого характера $\chi\in X_+\setminus\{1\}$ имеем $\theta\bot\chi E$ (так как $\chi E\subset\chi_1 E$), а потому $\int_G\theta\overline{ f\chi}\,dx=0.$ Значит,  $\widehat{f\overline{\theta}}=0,$ откуда $f=0$ п.в. Таким образом, $E=\theta H^2(G).$ Единственность  $\theta$ с точностью до множителя из $\mathbb{T}$ доказывается точно так же, как и равенство $\Theta =c\Theta_1$ в  доказательстве теоремы 3.6.$\Box$

\textbf{Замечание 5.3.} Как известно, в случае $G=\mathbb{T}$  всякое нетривиальное инвариантное подпространство --- неприводящее (см., например, \cite{Nik1},  доказательство следствия 1.4.1). Если же $G\ne \mathbb{T}$ и $G$ содержит наименьший положительный элемент $\chi_1,$ то  всегда существует нетривиальное  инвариантное подпространство, не являющееся неприводящим. Действительно,  тогда множество $X_+\setminus X^i$  непусто \cite[следствие 1]{SbMath}, и в силу леммы 1.1  подпространство $E_1$ пространства $H^2(G)$ с ортонормированным базисом $X_+\setminus X^i$ инвариантно, но $\chi_1 E_1=E_1.$

\textbf{Следствие 5.4}.    \textit{Пусть  $X$ содержит наименьший положительный элемент и  $f\in H^2(G).$ Порожденное функцией $f$ замкнутое инвариантное подпространство $E_f$ является неприводящим, если и только если}
$$
\int_G\log|f(x)|\,dx>-\infty.
$$

Доказательство. Если функция $f$ удовлетворяет указанному неравенству, то пространство $E_f$ имеет вид $\theta H^2(G),$ где $\theta$ --- внутренняя функция (\cite[теорема 8.5.2]{Rud}) и значит  является неприводящим.

Пусть теперь $\int_G\log|f(x)|\,dx=-\infty.$ Тогда   $\int_Gf(x)\,dx=0$ в силу \cite[теорема 8.4.1]{Rud}, а потому $\int_Gg(x)\,dx=0$ при всех $g\in E_f.$ Следовательно, пространство $E_f$ не имеет вида $\theta H^2(G),$ где $\theta$ --- внутренняя функция, и в силу теоремы 5.2 не является неприводящим.$\Box$

\

\begin{center}
\S 6. Некоторые приложения
\end{center}

\

   По аналогии с классическим случаем (см.  \cite[раздел 6.6]{Pel}), следствие 4.2 и теорема 2.4  дают следующую информацию
   о равномерной аппроксимации символов рациональными функциями от $\chi_1.$ Для $\varphi\in L^\infty(G)$  положим
$$
d_n(\varphi,G):=\mathrm{dist}_{L^\infty}(\varphi, \mathcal{R}_n\circ\chi_1+H^\infty(G)).
$$

\textbf{Следствие 6.1.} \textit{Пусть $X$ содержит наименьший положительный элемент $\chi_1,$ $\varphi\in L^\infty(G).$
 Последовательность $d_n(\varphi,G)$ принадлежит $\ell^p,$ ($0<p<\infty$), если и только если $\varphi=\varphi_1\circ \chi_1,$
 где $\varphi_1\in B_p^{1/p}.$
}

Применим теперь полученные результаты к  дискретной версии операторов Ганкеля над группами, рассмотренной в \cite[пример 2]{MK}.

Пусть $\nu\in \ell^2(X).$ Из теоремы Планшереля легко следует, что оператор $\mathcal{G}_\nu:\ell^2(X_+)\rightarrow \ell^2(X_-)$, определяемый равенством
  $$
  \mathcal{G}_\nu f(\xi):=\sum\limits_{\chi\in X_+}\nu(\xi\chi^{-1})f(\chi),\quad  \xi\in X_-,
  $$
  унитарно эквивалентен  оператору Ганкеля $H_{\varphi}$   с символом $\varphi=
  \mathcal{F}^{-1}\nu\in  L^2(G)$  и обратно  (здесь  $\mathcal{F}^{-1}:L^2(X)\to L^2(G)$ --- обратное преобразование Фурье на группе $X;$ подробности см. в \cite[теорема 5]{MK}).  Попутно отметим, что имеются
    тесные связи между   операторами  $\mathcal{G}_\nu$ и операторами Винера-Хопфа $\mathcal{T}_\nu$ над группами,
     изучавшимися в \cite{Adukov}. В частности,
  $\mathcal{T}_\nu=C_\nu-\mathcal{G}_\nu,$ где $C_\nu$ --- оператор свертки с $\nu.$

  Если  $\nu\in \ell^1(X),$ то
   оператор  $H_{\varphi}$ ограничен, так как $\varphi\in C(G),$ а потому ограничен и оператор $\mathcal{G}_\nu.$
  Применительно к данным операторам теоремы 3.1, 2.3 и 2.4 выглядят следующим образом.

  \textbf{Следствие 6.2.} \textit{Пусть $X$ содержит наименьший положительный элемент $\chi_1$, $\nu\in \ell^1(X),$ $\psi:=P_-(\mathcal{F}^{-1}\nu).$}

  1) \textit{Оператор $\mathcal{G}_\nu$ имеет конечный ранг, если и только если
$\psi=R\circ \chi_1,$ где $R$  есть рациональная функция, все полюсы которой лежат в $\mathbb{D},$
 и при этом}
$$
\mathrm{rank }\mathcal{G}_\nu=\mathrm{deg} R.
$$

2) \textit{Оператор $\mathcal{G}_\nu$ компактен, если и только если
$\psi=\varphi_1\circ \chi_1,$ где $\varphi_1\in VMO(\mathbb{T}).$  }

3) \textit{Оператор $\mathcal{G}_\nu$ принадлежит классу  Шаттена--фон Неймана  $\mathbf{S}_p,$ $0<p<\infty,$ если и только если
$\psi=\varphi_1\circ \chi_1,$ где $\varphi_1\in B_p^{1/p}.$  }

(Операторы вида $\mathcal{G}_\nu$ рассматривались и для недискретных групп, см. \cite[определение 7, теорема 4]{MK}.)

 Следуя \cite{MD1}, рассмотрим еще одну версию дискретных операторов Ганкеля на группах. Ниже $1_{\{\xi\}}$ обозначает
  индикатор одноточечного подмножества $\{\xi\}\subset X.$

\textbf{ Определение 6.3} \cite{MD1}.  Оператор $\Gamma:\ell^2(X_+)\to \ell^2(X_+)$,
определенный первоначально на финитных  функциях  на $X_+,$
называется {\it ганкелевым (оператором Ганкеля) в} $\ell^2(X_+),$ если его матрица в базисе $\{1_{\{\chi\}}:\chi\in X_+\}$
 зависит только от произведения характеров,
т.~е. существует  такая функция $a_\Gamma$ на $X_+,$ что для всех
$\chi,\xi\in X_+$ выполняется равенство
$$
\langle\Gamma 1_{\{\chi\}},1_{\{\xi\}}\rangle=a_\Gamma(\chi\xi)
$$
(угловые скобки здесь обозначают скалярное произведение в $\ell^2(X)$).

Оператор $T$ в $\ell_2(X_+)$ будет ганкелевым тогда и только тогда,
 когда при всех  $\chi$ из $X_+$  выполняются коммутационные соотношения
 \[
 {\cal S}_\chi^* T=T{\cal S}_\chi
 \]
\cite[лемма 1]{MD1}, где
\[
 {\cal S}_\chi f(\xi)=
 \begin{cases}f(\chi^{-1}\xi), \mbox{ если } \chi^{-1}\xi\in X_+,\\ 0 ,\quad\quad\quad \mbox{ если } \chi^{-1}\xi\notin X_+
 \end{cases}
 \]
есть оператор в $\ell_2(X_+)$ сдвига на характер $\chi\in X_+,$ ${\cal S}_\chi^*f(\xi)=f(\chi\xi)$ --- сопряженный оператор.

При этом ганкелев оператор $\Gamma$ в
$\ell^2(X_+)$ ограничен тогда и только
тогда, когда существует  такая функция   $\psi\in  L^\infty (G),$ что $a_\Gamma$ есть сужение на $X_+$  преобразования Фурье $\widehat{\psi}$ функции $\psi$ \cite[лемма 6]{MD1}.

\textbf{Лемма 6.4.} \textit{Ограниченный оператор $\Gamma$ в
$\ell^2(X_+)$
унитарно эквивалентен   ганкелеву оператору $H_\varphi$ с символом $\varphi\in L^\infty (G)$ и обратно,
и при этом}
$$
P_-\varphi=\sum_{\xi\in X_-}a_\Gamma(\overline{\xi\chi_1})\xi.
$$

Доказательство. Утверждение об унитарной эквивалентности доказано в \cite[лемма 5]{MD1}, причем
в процессе доказательства показано, что $a_\Gamma(\chi)=\widehat{\varphi}((\chi\chi_1)^{-1})$ при
$\chi\in X_+.$  Следовательно,
$$
P_-\varphi=\sum_{\xi\in X_-}\langle P_-\varphi,\xi\rangle\xi=\sum_{\xi\in X_-}\widehat{\varphi}(\xi)\xi
=\sum_{\xi\in X_-}a_\Gamma((\xi\chi_1)^{-1})\xi
$$
(угловые скобки здесь обозначают скалярное произведение в $L^2(X)$).$\Box$

С учетом леммы 6.4, теоремы 3.1, 2.3 и 2.4  применительно к данным операторам  могут быть сформулированы следующим образом.

  \textbf{Следствие 6.5.} \textit{Пусть $X$ содержит наименьший положительный элемент $\chi_1,$ $h:=\sum_{\xi\in X_-}a_\Gamma(\overline{\xi\chi_1})\xi$}

  1) \textit{Оператор $\Gamma$ имеет конечный ранг, если и только если функция
$h$ имеет вид $R\circ \chi_1,
$
 где $R$  есть рациональная функция, все полюсы которой лежат в $\mathbb{D},$ и при этом}
$$
\mathrm{rank }\Gamma=\mathrm{deg} R.
$$

2) \textit{Оператор $\Gamma$ компактен, если и только если
$h$ имеет вид $\varphi_1\circ \chi_1,$ где $\varphi_1\in VMO(\mathbb{T}).$  }

3) \textit{Оператор $\Gamma$ принадлежит классу  Шаттена--фон Неймана  $\mathbf{S}_p,$ $0<p<\infty,$ если и только если
  $h$ имеет вид $\varphi_1\circ \chi_1,$
 где $\varphi_1\in B_p^{1/p}.$  }

\end{document}